\documentclass[
12pt]{amsart}
\usepackage{amsmath,amssymb,amsthm,url}
\usepackage{graphicx}
\usepackage{epic,eepic}
\usepackage[russian]{babel}
 \def\Z{{\bf Z}}   \def\Z{\mathbb{Z}}
\def\ord{\mathop{\fam0 ord}}

\def\id{\mathop{\fam0 id}}
\long\def\comment#1\endcomment{}

\newtheoremstyle{mydefinition}% name
  {3pt}%      Space above
  {3pt}%      Space below
  {\normalfont}%         Body font
  {\parindent}%         Indent amount (empty = no indent, \parindent = para indent)
  {\bfseries}% Thm head font
  {.}%        Punctuation after thm head
  { }%     Space after thm head: " " = normal interword space;
  {}%         Thm head spec (can be left empty, meaning `normal')

\theoremstyle{mydefinition}
\newtheorem{pr}{}[]
\begin{document}

\centerline{{\bf WHEN ANY GROUP OF $N$ ELEMENTS IS CYCLIC?}
\footnote{%See update version on www.arxiv.org.
We would like to acknowledge M. Vyalyi, P. Kozhevnikov and K. Kohas for useful discussions.}
}
\smallskip
\centerline{\bf V. Bragin, Ant. Klyachko and A. Skopenkov
\footnote{Homepage: \url{www.mccme.ru/~skopenko}. Supported by Simons-IUM Fellowship}
}

\bigskip

\small
We give a simple proof of the well-known fact: any group of $n$ elements is cyclic if and only if
$n$ and $\phi(n)$ are coprime.
This note is accessible for students familiar with permutations and basic number theory.
No knowledge of group theory is required; a few necessary notions are introduced in the course of the proof.
%The note could also be an interesting easy reading for mature mathematicians.

\normalsize

\bigskip
{\bf Introduction}

We call a {\it group} a nonempty family $G$ of transformations (i.e. permutations or rearrangements)
of some set, which family is closed with respect to composition and taking inverse transformation
(i.e. if $f,g\in G$, then $f\circ g\in G$ and $f^{-1}\in G$).
Common term: transformation group.
Cf. [A, comment to problem 5].

If a finite group $G$ contains an element $g$ such that $G$ consists of all powers of $g$
(i.e. $G=\{g,g^2,\dots,g^n,\dots\}$, then group $G$ is called {\it cyclic\/}.

We give a simple proof of the following well-known fact.

\smallskip
{\bf Theorem.}
{\it Any group consisting of $n$ elements is cyclic if and only if $n$ and $\phi(n)$ are coprime.}

\smallskip
Here $\phi(n)$ is the number of positive integers not exceeding
$n$ and coprime to $n$ (the Euler function).
Note that $n$ and $\phi(n)$ are coprime if and only if
in the prime decomposition $n=p_1\dots p_k$

(*) all $p_i$ are different and

(**) $p_i$ does not divide $p_j-1$ for any $i$ and $j$.

%Our proof is accessible for students familiar with permutations and basic number theory.
%No knowledge of abstract group theory is required;
%In particular,
Although we did not find such a proof in the literature, we do not claim any novelty.
\footnote{Although we do not use the Sylow theorems, our argument in the second case below is similar to their proof.
We do not use the notion of a quotient group, as opposed to more traditional proofs (see, e.g., [B]).}
One can see from [BKKSS] (or from the complete Russian version of this text) how to invent this proof.

%%%%%%%%%%%%%%%%%%%%%%%%%%%%%%%%%%%%%%%%
\bigskip
{\bf Proof of the ``only if" part.}

If $g$ and $h$ are transformations of disjoint sets $M$ and $N$, then we can regard them as transformations of
$M\sqcup N$ and hence take their composition.

A cycle $(a_1,a_2,...,a_n)$ is the transformation of a set containing $a_1,a_2,\dots,a_n$ that carries
$a_n$ to $a_1$ and $a_i$ to $a_{i+1}$ for each $i<n$, whereas it carries every other element to itself.

If condition (*) above is violated, e.g., $p_1=p_2=p$, then the following group
consists of $n$ elements and is not cyclic:
$$\left\{(1,2,\dots,p)^i\circ(p+1,p+2,\dots,p+\frac n p)^j\ |\ i=1,\dots,p,\ j=1,\dots,\frac np\right\}.$$
\quad
Assume that condition (**) above is violated, e.g., $p_1$ divides $p_2-1$.
Denote by $\Z_k$ the set of residues modulo $k$ with the summation and the multiplication operations.
Then from the primitive root theorem it follows that there is $a\in\Z_{p_2}$ for which the powers $a,a^2,\dots,a^{p_1}=1$ are different.
Denote by $G_{p_1,p_2}$ the group of transformations
$f_{k,l}:\Z_{p_2}^2\to \Z_{p_2}^2$ defined by the formula
$f_{k,l}(x,y):=(a^kx,lx+y)$ for $k\in\Z_{p_1}$ and $l\in\Z_{p_2}$.
\footnote{In more advanced notation
$G_{p_1,p_2}:=
\left\{ \begin{pmatrix} a^k&l\cr 0&1\cr
\end{pmatrix}\in\Z_{p_2}^{2\times2}
\;\Big|\;
k\in\Z_{p_1},\ l\in\Z_{p_2}\;
\right\}.
$}
Then the following group is not cyclic (it is even nonabelian):
$$\left\{f\circ(1,2,\dots,\frac n{p_1p_2})^j\ |\ f\in G_{p_1,p_2},\ j=1,2,\dots,\frac n{p_1p_2}\right\}.\quad QED$$

%As opposed to $\Z/k$, the object $\Z_k$ is the {\it set} having summation and multiplication operations known from number theory.
%\newline
%Denote $\Z/k:=\{(12\dots k)^j\ |\ j=1,2,\dots,k\}$.
%If $g$ and $h$ are transformations of sets $M$ and $N$, respectively, then transformations
%$g\times h$ of the Cartesian product $M\times N$ is defined by $(g\times h)(x,y):=(g(x),h(y))$.
%Denote $G\oplus H=\{g\times h\ |\ g\in G,\ h\in H\}$.

%%%%%%%%%%%%%%%%%%%%%%%%%%%%%%%%%%%%%%%%
\bigskip
{\bf Proof of the ``if" part.}

Denote by $|X|$ the number of elements in a set $X$.
Denote given group by $G$.

We use the induction on the number of prime factors of $|G|$.
If $|G|$ is a prime, then the ``if" part is implied by the following Lagrange Theorem.

The {\bf order} $\ord a$ of an element $a$ of a group with the identity
element $e$ is the minimal positive integer $n$ such that $a^n=e$.
If the group is finite, it is clear that such $n$ exists.

\smallskip
{\bf Lagrange Theorem (particular case).}
{\it The number of elements of any finite group is divisible by the order of any its element.}

\smallskip
{\it Proof.} Denote the group by $G$.
For each $x\in G$ consider the set
$\{x,xf,xf^2,\dots,xf^{\ord f-1}\}$.
 By the definition of order these elements are different.
Therefore this set contains $\ord f$ elements.
If $xf^k=yf^l$, then $y=xf^{k-l}$.
Therefore for different $x$ these sets either coincide or are disjoint.
Thus $|G|$ is divisible by $\ord f$.
QED

\def\gp#1{\left\langle#1\right\rangle}

\smallskip
Now suppose that the number of prime factors in $|G|$ is greater than one.
We need the following general version of the Lagrange Theorem.

A {\bf subgroup} of a group $G$ is a subset of $G$ that is itself a group.

\smallskip
{\bf Lagrange Theorem.}
{\it The number of elements of any finite group is divisible by the number of elements of any subgroup.}

\smallskip
{\it Proof.} Denote the group by $G$ and the subgroup by $\{h_1,h_2,\dots,h_m\}$.
For each $x\in G$ consider the set $\{xh_1,xh_2,\dots,xh_m\}$.
 This set contains $m$ elements.
If $xh_k=yh_l$, then $y=xh_kh_l^{-1}$.
Therefore for different $x$ these sets either coincide or are disjoint.
Thus $|G|$ is divisible by $m$.
QED

\smallskip
A {\bf maximal subgroup} of a group is a maximal by inclusion subgroup not coinciding with $G$ and containing more than one element.
By the induction hypothesis and the Lagrange Theorem, {\it each maximal subgroup is cyclic.}

For an element $f$ of  a group $G$ let $\gp f$ be the set of all powers of $f$ (including zero and negative ones).
The element $f$ is called {\bf generating} for the (cyclic) subgroup  $\gp f$.

Suppose to the contrary that the group $G$ is noncyclic.
Then each element is contained in a maximal subgroup.

Elements $f,g$ of a group $G$ are {\bf conjugate} in $G$ if $g=b^{-1}fb$ for some $b\in G$.

\smallskip
{\bf First case:} {\it generator $f$ of some maximal subgroup is conjugate only to (some of) its powers.}
Take $h\in G\setminus\gp f$.
Let $q$ be the minimal positive integer $n$ such that $h^n\in \gp f$.
Such a $q$ exists because $h^{\ord h}\in \gp f$.

\smallskip
{\it Proof that $|G|$ is divisible by $q$.}
Let $\ord h=qt+r$ be division of $\ord h$ on $q$ with remainder $r$.
Then $h^r=h^{\ord h-qt}\in \gp f$ and $0\le r<q$.
Hence $r=0$ by the minimality of $q$.
So $\ord h$ is divisible by $q$.
Therefore by Lagrange Theorem $|G|$ is divisible by $q$.
QED

\smallskip
{\it Proof that $fh=hf$.} Since $h^{\ord h}\in \gp f$, using division with a remainder we obtain that
By the condition of the first case $h^{-1}fh=f^k$ for some $k\in\Z$.
The inclusion $h^q\in\gp f$ implies $f=h^{-q}fh^q=f^{k^q}$
(here the last equality holds for each $q$ and is proved by induction on $q$).
Therefore $k^q\equiv1\mod\ord f$.
Hence $k$ and $\ord f$ are coprime.
By conditions (*),(**) and the Lagrange Theorem $|G|$ and $\varphi(\ord f)$ are coprime.
Since $|G|$ is divisible by $q$, numbers $q$ and $\varphi(\ord f)$ are coprime.
So there are integers $x$ and $y$ such that $qx+\varphi(\ord f)y=1$.
Thus $k\equiv k^{qx+\varphi(\ord f)y}\equiv1\mod\ord f$.
Hence $fh=hf$.
QED

%Hence $\varphi(\ord f)$ is divisible by $q$. Hence $q=1$ and so $fh=hf$.

\smallskip
{\it Completion of the argument for the first case.}
Since $fh=hf$, the group $G$ contains a subgroup
$$\{f^ih^j\ |\ 1\le i\le \ord f,\ 1\le j\le q\}$$ of $q\ord f$ elements.
Hence by condition (*) and the Lagrange Theorem $\ord f$ is coprime to $q$.
Since $(fh)^j=f^jh^j$ for each $j$, we obtain that $\ord(fh)$ is divisible both by $q$ and by $\ord f$.
Hence $\ord(fh)=q\ord f$.
Thus $\ord(fh)=q\ord f$.
Since the subgroup $\gp f$ is maximal, we have $\gp{fh}=G$.
Thus $G$ is cyclic.
Contradiction. QED

\smallskip
{\bf Second case:} {\it generator of any maximal subgroup is conjugate not only to its powers.}

The {\bf product of subsets $X$ and $Y$} of a group $G$ is the
set of all products $xy$, where $x\in X$ and $y\in Y$.
If one of these subsets consists of only one element,
e.g., $Y=\{y\}$, then we write $Xy$ instead of $X\{y\}$.

\smallskip
{\it Assertion 1. Any maximal subgroup $F$ contains the} {\bf center}
$$
Z=Z(G):=\{a\in G\ :\ ga=ag\text{ for any }g\in G\},
$$
{\it i.e., the set of elements commuting with each element of the group.}

\smallskip
{\it Proof.} Otherwise $FZ$ is a larger commutative subgroup.
By the maximality of $F$ we have $FZ=G$.
Hence $G$ is commutative.
This contradicts to the assumption of the second case. QED

\smallskip
{\it Assertion 2. The intersection of two maximal subgroups equals the center.}

\smallskip
{\it Proof.}
A nontrivial element of the intersection commutes with all elements of both subgroups.
Hence it commutes with any product of several multiples, each multiple being an element of one of our subgroups.
The set of such products is a subgroup.
By the maximality of our subgroups this subgroup coincides with the entire group.
Therefore the intersection is contained in the center.

Assertion 1 implies the converse inclusion.  QED

\smallskip
{\it Conclusion of the proof of the second case: calculations.}
Recall that any element of $G$ is contained in certain maximal subgroup.
By Assertion 2 for any {\it non-central} element such a subgroup is unique.
So the group is split into the center and disjoint union of complements
of maximal subgroups to the center.
The number of elements in such complements are the same for conjugate subgroups.
Denote by $\widehat F$ the number  of non-central elements of $G$ in the union of subgroups conjugate
to given maximal subgroup $F$.
Let $F_1,\dots,F_s$ be a maximal family of pairwise non-conjugate maximal subgroups.
Then
$$|G|=|Z|+\sum_{i=1}^s \widehat{F_i}.$$
By the left inequality in the following Assertion 3 the number of summands is at most one;
by the right inequality one summand is also impossible.
QED

\smallskip
{\it Assertion 3. $|G|/2\le \widehat F<|G|-|Z|$.}

\smallskip
{\it Proof.}
A subgroup conjugate to a maximal subgroup is also maximal.
(Indeed, if $g^{-1}Fg\subset F'\subset G$, then $F\subset gF'g^{-1}\subset G$.)

Consider the set
$$
N(F):=\{a\in G\ :\ Fa=aF\}.
$$
Then the number of different subgroups conjugate to $F$ (including $F$) is $|G|/|N(F)|$.

(Indeed, the conjugation by each element of $G$ takes $F$ to a conjugate subgroup.
If the conjugation by two different elements $u$ and $v$ takes the $F$ to the same
subgroup, i.e.,  $u^{-1}Fu=v^{-1}Fv$, then $Fuv^{-1}=uv^{-1}F$.
This means that $uv^{-1}\in N(F)$ or, equivalently, $u\in N(F)v$.
Conversely, the condition $u\in N(F)v$ implies $u^{-1}Fu=v^{-1}Fv$.
Clearly, $|N(F)v|=|N(F)|$.
Therefore the number of elements of $G$ conjugation by which takes $F$ to a given
subgroup equals $|N(F)|$.
Therefore the number of different subgroups conjugate to $F$ is precisely $|N(F)|$ times less than $|G|$.)

We have $N(F)=F$.

(Indeed, it is easy to verify that $N(F)$ is a subgroup.
By the assumption of the second case $N(F)\ne G$.
Since $N(F)\supset F$, the maximality implies that $N(F)=F$.)

The three statements just proved imply that
$\widehat F=(|F|-|Z|)\dfrac{|G|}{|F|}=|G|\left(1-\dfrac{|Z|}{|F|}\right)$.

Since $|G|>|F|$, we have $\widehat F<|G|-|Z|$.

By the assumption of the second case, $Z\ne F$.
By Assertion 1 the center is a subgroup of $F$.
Hence  by the Lagrange theorem $|Z|$ divides $|F|$.
Therefore $\widehat F\ge|G|/2$.
QED

\bigskip
\centerline{\bf References}

[A] V. I. Arnold, Ordinary Differential Equations, The MIT Press (1978), ISBN 0-262-51018-9.

[B] Ken Brown, Mathematics 4340, When are all groups of order $n$ cyclic?
Cornell University, March 2009,
\url{http://www.cornell.edu/~kbrown/4340/cyclic_only_orders.pdf}

[BKKSS] D. Baranov, A. Klyachko, K. Kohas, A. Skopenkov and M. Skopenkov,
When are all groups of order $n$ cyclic?
\url{http://olympiads.mccme.ru/lktg/2011/6/index.htm}

\newpage
\centerline{\uppercase{\bf Когда любая группа из $N$ элементов циклическая?}
\footnote{
%Обновляемая версия поддерживается на http://arxiv.org/abs/1108.5406
Благодарим М. Вялого, О. Иванова, П. Кожевникова, К. Кохася, А. Сгибнева, Б. Френкина и А. Шеня и за полезные замечания.}
}
\bigskip
\centerline{\bf В. Брагин, Ант. Клячко и А. Скопенков
\footnote{Поддержан грантом фонда Саймонса.
%Независимый Московский Университет.
Инфо: \url{www.mccme.ru/~skopenko}}
 }

\bigskip
\small
{\bf Аннотация.}
{\it Группой} называется непустое семейство $G$ перестановок некоторого множества, замкнутое относительно
композиции и взятия обратной перестановки.
Приводится простое доказательство ответа на следующий вопрос:
{\it для каких $n$ в любой группе из $n$ элементов найдется перестановка $g\in G$, для которой
$G= \{g,g^2,\dots,g^n\}$?}
%, тогда и только тогда, когда $n$ взаимно просто с $\varphi(n)$.}
Для понимания доказательства необходимо знание основ теории чисел (включая теорему Ферма-Эйлера).
Знаний по теории групп не требуется: небольшое количество понятий, необходимых для доказательства, приведены.
\normalsize

\bigskip
\hfill{\it It startled the well informed by being a new and fantastic}

\hfill{\it idea  they had never encountered. It startled the ignorant by being}

\hfill{\it an old and familiar idea they never thought to have seen revived.}

\hfill{\it G. K. Chesterton. The Man Who Knew Too Much.}

%\hfill{\it Like the witch, he liked to answer a question with a question; but the}
%\hfill{\it answers to Rose's questions were always something she'd always known, while}
%\hfill{\it ...The answers to his questions were things she had never imagined}
%\hfill{\it and found startling, unwelcome, even painful, altering her beliefs.}
%\hfill{\it U. K. Le Guin, Dragonfly.}

\tableofcontents

\section{Введение: зачем, для кого и как устроена эта заметка}

Мы хотели бы привлечь внимание к теории групп широкого круга людей, включая учителей, руководителей кружков и школьников, серьезно интересующихся математикой и программированием.
В этой теории есть доступные и интересные им результаты-жемчужины.
Формулировки таких результатов кратки и используют лишь простейшие определения;
доказательства красивы и похожи на решения сложных олимпиадных задач.
К сожалению, в б\'ольшей части существующей литературы эти жемчужины погребены под огромным количеством
немотивированного материала.
%, что делает их неинтересными и недоступными.
На примере исследования просто формулируемого вопроса мы покажем, как {\it появляются} некоторые основные понятия теории групп.
Основные идеи представлены на <<олимпиадных>> примерах: на простейших частных случаях,  свободных от технических деталей.
\footnote{Это не только делает материал более доступным, но помогает тем,  кто привык к  абстрактному изложению, развить математический вкус.
Благодаря этому они смогут разумно выбирать проблемы для исследования и ясно излагать собственные открытия, не скрывая ошибок (или известности полученного результата) за чрезмерным формализмом.
К сожалению, такое (непреднамеренное) сокрытие иногда происходит с математиками, воспитанными на чрезмерно формальных курсах.}
%Происходило и с третьим автором заметки; к счастью, почти все его ошибки исправлялись {\it перед} публикациями.}

Эта заметка предназначена для того, кому понятна и интересна формулировка вопроса, сформулированного в аннотации.
Мы приведем простое доказательство теоремы, отвечающей на этот вопрос.
Оно не претендует на новизну, хотя мы не нашли такого доказательства в литературе.

Заметка может быть интересна читателю, не знакомому с основами абстрактной теории групп, но изучавшему перестановки и основы теории чисел.
Например, по книгам [Al, KS, V] и [GDI, п. 1.7, 1.8 и 1.9].
Она может быть интересна и читателю, знакомому с этими основами, ибо ответ на сформулированный вопрос нетривиален.
Такому читателю может оказаться достаточным прочитать \S3.

Эта заметка не должна быть единственным и даже первым шагом в теорию групп.
В миникурс <<Рождение понятия группы>> можно включить [GDI, п. 1.7, 1.8 и 1.9],
%{\it (эти пункты предполагается включить в новое издание [Z])},
[Z, глава 6], [Z, глава 3, разделы <<Малая теорема Ферма>>, <<Квадратичные вычеты>> и <<Первообразные корни>>], [S08], [S15],
%(там приведено простое доказательство {\it основного} результата книги [Al]),
данную заметку и другой материал.

% (и имеющему склонность к задачам классификации).
Для понимания доказательства необходим опыт работы с перестановками  и числами (включая теорему Ферма-Эйлера).
% (на элементарном языке).
\footnote{В частности, наше доказательство не привлекает явно понятия факторгруппы,
в отличие от более традиционных доказательств, см., например, [B].
Также, мы не используем теорем Силова, хотя наш разбор второго случая похож на их доказательство.}
Знаний по теории групп и опыта работы с определением абстрактной группы не требуется.
Небольшое количество необходимых понятий вводятся (и могут быть освоены читателем) в процессе доказательства.

Конечно, читателю, не знакомому с основами теории групп, нужно будет потрудиться, чтобы самостоятельно доказывать тривиальные факты про эти понятия.
Такие упражнения --- важная часть изучения этих понятий.
Выполнять их интереснее ради красивых результатов, формулировки которых ясны и доступны неспециалисту
(в частности, не используют этих понятий), но в доказательствах которых эти понятия возникают.
%, чем ради сдачи зачета или изучения немотивированной теории.
%Т.е. в качестве части {\it пути к} пониманию теории, а не в качестве {\it награды за}
А не в процессе долгого немотивированного изучения
%этой
теории.
Эта заметка будет особенно интересна читателю, предпочитающему  изучить доказательство красивого результата на несколько страниц, самостоятельно разбираясь в деталях, чем прочитать сотню страниц более легкого материала,
%большая часть которых
не мотивированных таким результатом.

Опыт работы с абстрактными группами как раз появится при изучении данной заметки, хотя в ней формально не используется этого понятия.
Читатель увидит, что помимо всех рассматриваемых объектов есть еще одно множество (на котором действуют перестановки).
Странным образом оно так никогда и не выходит из тени.
В результате естественно возникает общее понятие группы.
%к концу доказательства для большинства читателей общее понятие группы уже естественно вырисовывается.
Итак, {\it данная заметка посвящена мотивировке важного общего понятия группы}
(здесь оно не используется, но его и основы соответствующей теории можно найти в [Al, K, KS]).
\footnote{{\it Добавление от А. Скопенкова.}
О необходимости мотивировок говорили классики математики [K], [P1, Глава 2 `Математические определения и преподавание'], [P2, стр. 455-475], [A, стр. 49, комментарий к задаче 5];
см. также [Z, Философски-методическое отступление].
%, увы, обычной ситуации: рассказчик (возможно, убеждается в полезности изученных им понятий и потому забывает,
%что полезность нужно обосновывать --- желательно конкретными результатамии
%и при рассказе данного материала другим, и в собственных научных работах.
Процитируем В.И. Арнольда:
%Алгебраисты
`...Обычно определяют группу как множество с двумя операциями, удовлетворяющими набору аксиом вроде $f(gh) = (fg)h$.
Эти аксиомы автоматически выполняются для групп преобразований.
В действительности эти аксиомы означают просто, что группа образована из некоторой группы преобразований забыванием преобразуемого множества.
Такие аксиомы, наряду с другими немотивированными определениями, служат математикам главным образом для того, чтобы затруднить непосвященным овладение своей наукой и тем самым повысить ее авторитет.'}
Поэтому усилия по изучению заметки будут сполна вознаграждены тем, что вслед за великими математиками в~процессе изучения конкретной задачи читатель увидит, как естественно возникает важное понятие.
Надеюсь, это поможет ему совершить собственные настолько же полезные открытия (не обязательно в~математике)!

% Именно с таких жемчужин полезно начинать изучение абстрактной теории групп,
%на примере их доказательства показывая, как появляются ее основные понятия.

Параграфы 2 и 3 формально независимы друг от друга.

%\newpage
\section{Как придумать}
 %\footnote{Это пункт может быть пропущен читателем, знакомым с основами теории групп.}

\subsection{Постановка задачи}

{\it Общие замечания к формулировкам задач.} Eсли условие задачи является
утверждением, то в задаче требуется это утверждение доказать.
Если некоторая задача не получается, то читайте дальше --- соседние задачи могут оказаться подсказками.
%Из предыдущего текста используется только пункт `Примеры' и только во введении.

\begin{pr}\label{open}
Дано семейство $G$ из 11 перестановок некоторого множества, замкнутое относительно композиции и взятия обратной  перестановки (т.е. если $f,g\in G$, то $f\circ g\in G$ и $f^{-1}\in G$).
Тогда найдется перестановка $g\in G$, для которой $G= \{g,g^2,\dots,g^{11}\}$.
\end{pr}

\begin{pr}\label{pr-cycle}
Верен ли аналог предыдущего утверждения для аналогичного семейства $G$ из $n$ перестановок при $n=$

(7) 2,3,4,5,6,7 (ответ может быть разным для разных значений $n$); \quad

(8) 8; \quad (9) 9; \quad (10) 10; \quad (12) 12;
\quad (15) 15; \quad (21) 21; \quad (1001) 1001?
\end{pr}

{\bf Основной Вопрос.} {\it Дано семейство $G$ из $n$ перестановок некоторого множества, замкнутое относительно композиции и взятия обратной перестановки.
Для каких $n$ обязательно найдется перестановка $g\in G$, для которой $G= \{g,g^2,\dots,g^n\}$?}

\smallskip
{\it Группой} называется непустое семейство $G$ преобразований (т.е.  перестановок) некоторого множества, замкнутое относительно композиции и взятия обратного преобразования (т.е. если $f,g\in G$, то $f\circ g\in G$ и $f^{-1}\in G$).
\footnote{Это определение равносильно обычному ввиду теоремы Кэли.
Группы, определенные как здесь, обычно называют {\it группами преобразований}.}

Если в конечной группе $G$ найдется перестановка $g$, из всех возможных степеней которой состоит $G$
(т.е. $G=\{g,g^2,\dots,g^n,\dots\}$), то эта группа называется {\it циклической\/}.
Примеры циклических и нециклических групп Вы привели при решении задачи \ref{pr-cycle}.

\begin{pr}\label{pr-anyn} Для любого $n$ имеется циклическая группа из $n$ элементов.
\end{pr}

\subsection{План \S2}

Пункт \ref{s:motiv} нужен только для мотивировки, задачи \ref{pr-constr}.c и п. 3.2.
Пункт \ref{s:exa} подводит читателя к построению примеров, необходимых для ответа на основной вопрос, т.е., к доказательству достаточности критерия, который мы нащупаем.
(Немногие задачи в пунктах \ref{s:lagr}-\ref{s:zentr}, использующие пункт \ref{s:exa}, можно пропустить, ибо эти задачи не касаются основного вопроса.)
А вот пункты \ref{s:lagr}-\ref{s:zentr} подводят читателя к необходимым условиям из этого  критерия.  
Пункты \ref{s:exa}-\ref{s:zentr} интересны и сам по себе, независимо от основного вопроса. 

Пункты \ref{s:motiv}, \ref{s:exa} и \ref{s:lagr}-\ref{s:zentr} формально независимы друг от друга. 
В пунктах \ref{s:exa}-\ref{s:zentr} используется следующее определение и обозначение.  

{\it Циклом} $(a_1,a_2,...,a_n)$ называется перестановка множества, содержащего элементы
%\linebreak
$a_1,a_2,\dots,a_n$,
которая переводит $a_n$ в $a_1$ и $a_i$ в $a_{i+1}$ для любого $i<n$, а каждый из остальных элементов переводит в себя.

\subsection{Почему вопрос интересен}\label{s:motiv}

\smallskip
{\bf Теорема о первообразном корне.}
\footnote{Операция умножения на множестве ненулевых вычетов по простому модулю имеет общее обобщение с операцией композиции перестановок. Но для доказательства результата заметки не нужно понимать этого.}
{\it Число $g$, для которого $g^1,g^2,g^3,\dots,g^{\varphi(n)}$ есть
все вычеты по модулю $n$, взаимно простые с $n$, существует тогда и только тогда, когда $n$ есть либо 2, либо 4, либо степень нечетного простого, либо удвоенная степень нечетного простого.}

%остатки от деления на $m$ чисел $g^1,g^2,g^3,\dots,g^{\varphi(m)}$ различны.}

\smallskip
Поставленный вопрос касается аналога теоремы о первообразном корне для перестановок.
Впрочем, ситуация для перестановок не совсем аналогична: существуют циклическая и нециклическая группы с одинаковым числом элементов.

В процессе изучения этого вопроса мы в основном будем работать преобразованиями, забыв про множество, на котором они действуют.
Такая ситуация часто встречается в математике и мотивирует следующее определение.
 Группы преобразований $G$ и $H$ множеств $M$ и $N$ называются {\it изоморфными}, если существует биекция
$\varphi:G\to H$, для которой $\varphi(g_1)\circ\varphi(g_2)=\varphi(g_1\circ g_2)$ при любых $g_1,g_2\in G$.
Важная тема в математике --- классификация групп с точностью до изоморфизма.
Ясно, что для любого $n$ имеется циклическая группа из $n$ элементов, ровно одна с точностью до изоморфизма.
Поэтому вопрос о {\it цикличности любой группы из $n$ элементов} равносилен вопросу о
{\it единственности группы из $n$ элементов с точностью до изоморфизма}.

\subsection{Примеры групп}\label{s:exa}

(1) Группа $S_n$ {\it всех} перестановок $n$-элементного множества.

(2) Группа  перестановок $\{\id=(1)(2)(3)(4),(13)(24),(1234),(1432)\}$
множества из четырех элементов.

(3) Группа перестановок $\{(1,2,3,4,5,6,7,8)^k\}$, $k=1,2,\dots,8$,
%восьми
8-элементного множества.

(4) Группа перестановок $\{(1,2,3)^k(4,5,6,7,8)^l\}$, $k=1,2,3$, $l=1,2,3,4,5$,
 8-элементного множества.

%\smallskip
\noindent
\centerline{\begin{tabular}{@{}c@{}c@{}}
\includegraphics{skopenkov-003.mps}\qquad\includegraphics{skopenkov-004.mps}\\
Рисунок: движения квадрата и куба\\
\end{tabular}}
%\smallskip

(5) Рассмотрим квадрат на плоскости и все движения плоскости, переводящие его в себя.
Это тождественное преобразование, 3 поворота и 4 симметрии.
Всего 8 преобразований.
Возьмем группу из 8 перестановок множества вершин квадрата,
происходящих при применении перечисленных восьми преобразований плоскости.

(6) Рассмотрим куб в пространстве и все вращения пространства (включая тождественное), переводящие его в себя.

\quad (a) Возьмем группу из всех перестановок множества {\it вершин} куба,
происходящих при применении таких вращений.

\quad (b) Возьмем группу из всех перестановок множества {\it середин ребер} куба,
происходящих при применении таких вращений.

\smallskip
\noindent
\centerline{\begin{tabular}{@{}c@{}}
\includegraphics{skopenkov-002.mps}\\
Рисунок: граф $K_{3,3}$\\
\end{tabular}}
\smallskip

\smallskip
\noindent
\centerline{\begin{tabular}{@{}c@{}}
\includegraphics{skopenkov-001.mps}\\
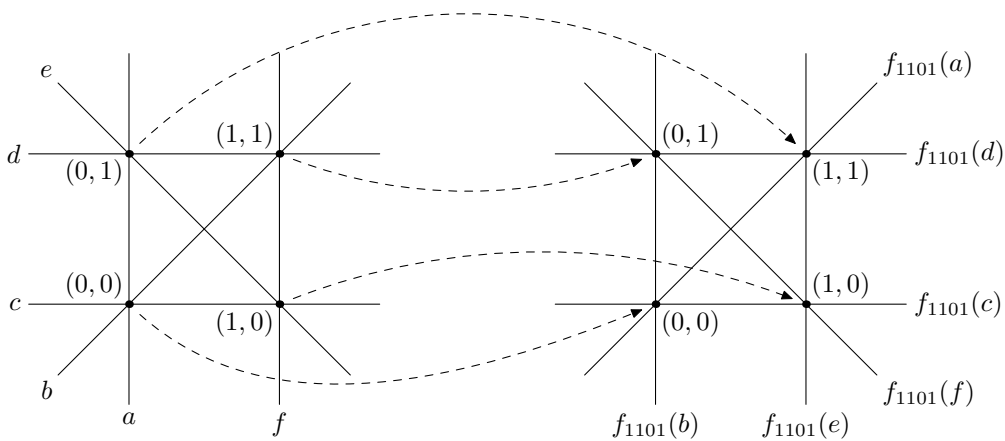
Рисунок: линейное преобразование $f_{1101}:\Z_2^2\to\Z_2^2$\\
\end{tabular}}
\smallskip

(7) Группа всех перестановок 6-элементного множества, являющихся изоморфизмами графа $K_{3,3}$.
({\it Изоморфизм} графа $G$ --- такая перестановка множества его вершин, что для любых двух вершин эти вершины  соединены ребром тогда и только тогда, когда их образы при перестановке соединены ребром.)

(8) Рассмотрим множество $\Z_2^2=\{(0,0),(0,1),(1,0),(1,1)\}$ упорядоченных пар вычетов по модулю 2.
Для любых четырех вычетов $a,b,c,d$ по модулю 2 рассмотрим отображение $f_{abcd}:\Z_2^2\to\Z_2^2$, 
заданное формулой $f_{abcd}(x,y)=(ax+by,cx+dy)$.
Среди всех таких отображений выберем взаимно-однозначные.
Они образуют группу. 
% преобразований.

\begin{pr}\label{pr-id}
(a) Любая группа содержит тождественное преобразование.
% (оно называется {\it единичным элементом}).

(b) Множество из примера 8 действительно является группой.
\end{pr}

\begin{pr}\label{pr-ex}
Какие из приведенных примеров групп являются циклическими?
\end{pr}

\begin{pr}\label{pr-constr}
(a) Cуществует нециклическая группа из 21 элемента.

(b) Cуществует нециклическая группа из 55 элементов.

(c) Если  $p$ и $q$ простые числа и $q-1$ делится на $p$, то существует
нециклическая группа из $pq$ элементов.
\end{pr}

%\newpage
\subsection{Докажем и применим теорему Лагранжа}\label{s:lagr}

\begin{pr}\label{pr-prime-cycle}
Если количество перестановок в группе является простым числом, то эта группа циклическая.
\end{pr}

{\it Порядком} $\ord a=\ord_Ga$ перестановки $a$ в группе $G$ называется наименьшее
целое положительное $n$, для которого $a^n=\id$ (если такое $n$ существует).

\def\gp#1{\left\langle#1\right\rangle}

\begin{pr}\label{pr-ordfh}
Пусть $G$ --- группа и $f,h\in G$.

(a) $f^n=\id$ тогда и только тогда, когда $n$ делится на $\ord f$.

(b) Число $\ord h$ делится на наименьшее $\ord_fh$ из целых положительных $n$, для которых $h^n\in \gp f:= \{f,f^2,f^3,\dots\}$.
\end{pr}

\begin{pr}\label{pr-finite-ord}
(a) Найдите порядок каждой перестановки в группе $S_4$.

(b) Любая перестановка конечной группы имеет (конечный) порядок.

(c) Если в конечной группе есть перестановка порядка 2, то число перестановок в группе четно.

(d) Если в конечной группе есть перестановка порядка 3, то число перестановок в группе делится на 3.

(e) {\it Теорема Лагранжа.} Число перестановок конечной группы делится на порядок любой ее перестановки.

(Это `некоммутативная версия' теоремы Ферма-Эйлера, а также подсказка к задаче \ref{pr-prime-cycle}.)
%(f) В любой группе из четного числа элементов есть элемент порядка 2.
\end{pr}

\begin{pr}\label{pr-contrprimer} (a) Если число $n$ четное составное,
то существует группа из $n$ перестановок, не являющаяся циклической.

(b) Если число $n$ делится на квадрат простого,
то существует группа из $n$ перестановок, не являющаяся циклической.
\end{pr}

Группа $G$ называется {\it коммутативной}, если $xy=yx$ для любых $x,y\in G$.

\begin{pr}\label{pr-comm=cyclic}

(b) Любая циклическая группа является коммутативной.

(c) Верно ли обратное?
\end{pr}

%{\it Теорема Ферма-Эйлера.} Для любого элемента $a$ конечной коммутативной группы $G$ с единичным элементом $e$ выполнено $a^{|G|}=e$.

\begin{pr}\label{pr-comm}
(a) Любая коммутативная группа из 10 перестановок является циклической.

(b) То же для 21 перестановки.

(c) То же для 1001 перестановки.

(d) Для каких $n$ любая коммутативная группа из $n$ перестановок является циклической?
\end{pr}

\begin{pr}\label{pr-10}
Может ли в коммутативной группе из 10 перестановок быть две различных перестановки порядка 2?

(Это подсказка к задаче \ref{pr-comm}.a.)
\end{pr}

\begin{pr}\label{pr-abel}*
Для каких $n$ любая группа из $n$ элементов является коммутативной?
(Решение этой непростой задачи лучше отложить до разрешения основного вопроса.)
\end{pr}

{\it Подгруппой} группы $G$ называется подмножество группы $G$, также являющееся группой.

\begin{pr}\label{pr-LagrTheorem}
{\it Теорема Лагранжа.} Число перестановок в конечной группе делится на число перестановок в любой ее подгруппе.

(Это подсказка к задаче \ref{pr-10}.)
\end{pr}

\begin{pr}\label{pr-orbits}
%(Это подсказка к задаче \ref{pr-sopr}b.)
Пусть $G$ --- группа из 15 элементов и $f,g\in G$ --- элементы порядка 5.

(a) Множества $\{f,f^2,f^3,f^4\}$ и $\{g,g^2,g^3,g^4\}$ пересекаются.

(b) Множества $\{f,f^2,f^3,f^4\}$ и $\{g,g^2,g^3,g^4\}$ совпадают.

(c) Один из элементов $f,g$ является степенью другого.

(d) В $G$ есть элемент порядка 3.
% (даже 10 штук НЕТ, 2).
\end{pr}

\begin{pr}\label{pr-Lagr} Пусть $G$ --- группа из 15 элементов  и $f,g\in G$ --- элементы порядка 3.

(a) Множества $\{f,f^2\}$ и $\{g,g^2\}$
либо не пересекаются, либо совпадают.

(b) Если $\{f,f^2\}\ne\{g,g^2\}$, то $fg\ne gf$.
\end{pr}

 %\newpage
\subsection{Сопряжение}

Цикличность любой группы из 15 перестановок вытекает либо из задач \ref{pr-orbits}.c, \ref{pr-fifteen}.b, \ref{pr-1001_2}.c и \ref{pr-cen}.c,
либо из задач \ref{pr-1001_2}.c и \ref{pr-alio}.d (второй способ сложнее, но его удобно переносить на общий случай).

Перестановки $f,g\in G$ называются {\it сопряженными} в группе $G$, если
$g=b^{-1}fb$ для некоторой перестановки $b\in G$.

\centerline{\begin{tabular}{@{}c@{}}
\includegraphics{skopenkov-005.mps}\\
Рисунок: перестановка типа $\left<1,2,3,4\right>$\\
\end{tabular}}

\begin{pr}\label{br-type}
(a) Перестановка $(n_1+...+n_k)$-элементного множества, являющаяся композицией непересекающихся циклов порядков $n_1,...,n_k$, называется перестановкой {\it типа} $\left<n_1,...,n_k\right>$.

Перестановки $f$ и $g$ сопряжены в $S_n$ тогда и только тогда они одного типа.

(b) Пусть $a$ и $x$ --- произвольные перестановки $n$-элементного множества. Тогда
$$xax^{-1}=\left(\begin{array}{cccc}
x(1) & x(2) & \ldots & x(n)\\
x(a(1))& x(a(2)) & \ldots & x(a(n))
\end{array}\right).$$
Иными словами, циклическое разложение перестановки $xax^{-1}$ получается из
циклического разложения перестановки $a$ заменой каждого элемента на его $x$-образ: если
\linebreak
$a=\Pi_{j=1}^q(i_{j1},i_{j2},\ldots,i_{js_j})$, то
$xax^{-1}=\Pi_{j=1}^q(x(i_{j1}),x(i_{j2}),\ldots,x(i_{js_j}))$.

(c) Вращения куба вокруг больших диагоналей сопряжены.
\end{pr}

\begin{pr}\label{pr-fifteen}
Пусть $G$ --- группа, $b,f\in G$ и $k$ --- целое положительное.

(a) $(b^{-1}fb)^k=b^{-1}f^kb$.

(b) Порядки сопряженных перестановок равны.

(c) Если $b^{-1}fb=f^m$, то $b^{-1}f^kb=f^{km}$ и $b^{-k}fb^k=f^{m^k}$
\end{pr}

\begin{pr}\label{pr-1001_2}
Пусть $G$ --- группа из 15 перестановок, $f\in G-\{e\}$ сопряжена только со
%некоторыми
своими степенями и $h\in G-\gp f$.

(a) $h^{-1}fh=f$.

(b) $\{f^ih^j\ |\ 1\le i\le \ord f,\ 1\le j\le q\}$ подгруппа в $G$.

(c) $G=\gp fh^{\ord f}$.
\end{pr}

Подгруппа группы $G$ называется {\it собственной}, если она не совпадает ни с $\{\id\}$, ни с
%всей
$G$.

\begin{pr}\label{pr-alio}
Пусть $G$ --- группа из 15 перестановок.
Пусть любая перестановка сопряжена не только со своими степенями.

(a) Любые две различные собственные подгруппы пересекаются только по тождественной перестановке.

(b) С собственной подгруппой из $s$ элементов сопряжено ровно $15/s$ различных подгрупп.

(c) Для числа $\widehat F$ нетождественных перестановок, сопряженных перестановкам данной собственной подгруппы, справедливы неравенства $7<\widehat F<14$.

(d) Получите противоречие.
 \end{pr}

Через $|X|$ обозначается число элементов в множестве $X$.

\begin{pr}\label{pr-cen}
Для группы $G$ и $f\in G$ обозначим через

$\bullet$ $c_G(f)$ число перестановок, сопряженных с $f$ (вместе с самой $f$).

$\bullet$ $z_G(f)$ число перестановок, коммутирующих с $f$ (вместе с самой $f$).

(a) Найдите $c_{S_3}(f)$ и $z_{S_3}(f)$ для каждого $f\in S_3$.

(b) $c_G(f)z_G(f)=|G|$.

(c) Если в группе из 15 перестановок нет перестановок порядка 5, то с каждой нетождественной перестановкой
сопряжено ровно 4 других перестановки.
\end{pr}

\begin{pr}\label{pr-normalizer}
(Это подсказка к \ref{pr-alio}.b.)
Для группы $G$ и $f\in G$ обозначим через
$n_G(f)$ число перестановок, сопряжение $f$ с которыми переводит $f$ в ее степень:
$$n_G(f):=|\{b\in G\ :\ fb=bf^k\text{ для некоторого }k\}|.$$
\quad
(a) Найдите $n_{S_3}(f)$ для каждого $f\in S_3$.

(b) $n_G(f)\ge\ord f$.

(c) $n_G(f)=\ord f$, если элемент $f$ сопряжен не только со своими степенями и $|G|=15$.

(d) Число различных подгрупп, сопряженных с $\gp f$, равно $|G|/n_G(f)$.
\end{pr}

\begin{pr}\label{pr-pq} Если число перестановок в группе есть произведение $pq$
простых чисел, $p<q$ и $q-1$ не делится на $p$, то эта группа циклическая.
\end{pr}

\begin{pr}\label{pr-case1}
Если $G$ --- группа из 1001 перестановки и $f\in G-\{e\}$ сопряжен только со
%некоторыми
своими степенями, то $G$ циклическая.
\end{pr}

%255, 455

%pr-fifteen
%(с) Любой элемент порядка 5 в группе из 15 элементов сопряжен только со
%%некоторыми
%своими степенями.

%pri-fifteen
%(c) Следуют из \ref{pr-orbits}.c.

\subsection{Максимальные подгруппы и центр}\label{s:zentr}

\begin{pr}\label{pr-main}
Пусть $G$ --- группа из 1001 элемента, не являющаяся циклической.

(a) Каждый элемент содержится в максимальной по включению подгруппе, не совпадающей со всей $G$.
Такие подгруппы будем сокращенно называть {\it максимальными подгруппами.}

(b) Каждая максимальная подгруппа является циклической.
\end{pr}

Максимальные подгруппы аналогичны собственным подгрупп, на рассмотрении которых
было основано доказательство для $n=15$.
Различные собственные подгруппы пересекались только по единичному элементу.
Доказать аналогичное утверждение для максимальные подгрупп при $n=1001$ не получается.
Как же устроено пересечение максимальных подгрупп?

\begin{pr}\label{pr-zCENTR}
Назовем {\it коммутативизатором} группы $G$ множество
$$Z=Z(G):=\{a\in G\ :\ ga=ag\text{ для любого }g\in G\}$$
тех элементов, которые коммутируют со всеми.
(Мы надеемся, что использование слова {\it коммутативизатор}
вместо общепринятого {\it центр} более удобно для начинающих.)

(a) Найдите $Z(S_n)$ для каждого $n=2,3,4,\dots$

(b) Коммутативизатор группы является подгруппой.
\end{pr}

\begin{pr}\label{pr-main2}
Пусть $G$ --- группа из 1001 элемента, не являющаяся циклической.
Пусть порождающий элемент любой максимальной подгруппы сопряжен не только со своими степенями.
%Напомним, что по утверждению задачи \ref{pr-main}b они циклические.

(a) Любая максимальная подгруппа содержит коммутативизатор.

(b) Пересечение двух максимальных подгрупп равно коммутативизатору.

(c) Число различных подгрупп, сопряженных с максимальной подгруппой из $s$ элементов, равно $1001/s$.

(d) Для числа $\widehat F$ элементов, сопряженных элементам максимальной подгруппы $F$ и не лежащих в коммутативизаторе, справедливы неравенства $500<\widehat F\le 1000-|Z|$.

(e) Получите противоречие.
\end{pr}

%(b) Выразите через $|F|$ и $|Z|$ число элементов, сопряженных элементам максимальной подгруппы $F$.

 %\newpage
\subsection{Указания, решения и ответы к некоторым задачам}

\smallskip
{\bf  \ref{open}.} Аналогично решению нижеприведенной задачи \ref{pr-cycle}.7 для $n=3$.

\smallskip
{\bf  \ref{pr-cycle}.} (7) для $n=3$. Пусть, напротив, имеется нециклическая группа $G$ из трех перестановок.
Обозначим через $a$ нетождественную перестановку в ней.
Если $a^2\ne e$, то перестановки $a,a^2,a^3$ различны и группа циклическая.
Если же $a^2=e$, то рассмотрим перестановку $b\in G$, отличную от $e$ и $a$.
Тогда перестановка $ab$ отлична от $e,a,b$.
(Действительно, $ab\ne a$ и $ab\ne b$ очевидно. Если $ab=e$, то $b=a^2b=a$ --- противоречие.)
Противоречие.

(10)  Рассмотрим правильный пятиугольник на плоскости.
Рассмотрим все движения плоскости, переводящие его в себя.
Это тождественное преобразование, 4 поворота и 5 симметрий.
Всего 10 преобразований.
Нужную группу образуют 10 перестановок множества вершин правильного пятиугольника,
происходящих при применении перечисленных десяти преобразований плоскости.
Эта группа нециклическая, поскольку для двух перестановок $s$ и $t$, <<пришедших из симметрий>>,
$st\ne ts$. А если бы $s=g^k$ и $t=g^l$ для некоторой перестановки $g$, то $st=ts$.

%oshibka pri kompiljacii, pdf normalnyj
%\centerline{\begin{tabular}{@{}c@{}c@{}}
%\includegraphics{fig5gon.pdf}\\[5pt]
%\label{fig-5gon}
%Рисунок: движения правильного 5-угольника\\
%\end{tabular}}
%\bigskip

%\begin{figure}[h]
%\input{fig5gon.tex}
%\caption{движения правильного 5-угольника} \label{fig-5gon} \end{figure}

%\smallskip
{\it Другое решение задачи  {\ref{pr-cycle}.10.}}
Оно более сложное, чем предыдущее, но зато может помочь Вам в решении задачи \ref{pr-constr}.
Рассмотрим перестановки $r$ и $s$ 5-элементного множества вершин правильного 5-угольника, <<сооответствующие>>  повороту на $2\pi/5$ и симметрии.
%; см. рисунок. %(рис.~\ref{fig-5gon}).
Тогда $r^5=e=s^2$ и $sr=r^{-1}s$.
Рассмотрим 10 перестановок $r^ks^l$, $k,l\in\Z$.
Из соотношения $sr=r^{-1}s$ можно получить, что это множество является группой.
(Именно в этом отличие приводимого решения от предыдущего --- мы получили замкнутость относительно композиции и взятия обратного не из геометрических соображений, а из комбинаторных.
Поэтому появилась возможность обобщать это доказательство на случаи, когда геометрической интерпретации не видно.)
Из этого же соотношения вытекает, что эта группа не является циклической.

\smallskip
{\bf  \ref{pr-cycle}.} Ответы: верен для $n\in\{2,3,5,7,15,1001\}$ и неверен иначе.

Доказательство верности для простых $n$ аналогично приведенному для $n=7$,
а для $n=15$ и $n=1001$ намечено в задачах далее.

Доказательство верности для четных $n>2$ аналогично приведенному для $n=10$,
а для $n=21$ намечено в задачах далее.

\smallskip
{\bf \ref{pr-anyn}.} Множество степеней цикла длины $n$.

\smallskip
{\bf \ref{pr-id}.}
\quad $f\in G\quad\Rightarrow\quad f^{-1}\in G \quad\Rightarrow\quad ff^{-1}=e\in G$.

\smallskip
{\bf  \ref{pr-ex}.} Ответ: группы из примеров (1) для $n=2$, (3) и (4) --- циклические, а из остальных примеров --- нет.

\smallskip
{\bf \ref{pr-constr}.} (a) Искомая группа является группой некоторых перестановок 49-элементного множества $\Z_7^2$.
Для описания группы представим его элементы  как пары $(x,y)$ вычетов по модулю 7.
Для любых целых неотрицательных $k,l$ определим преобразование
$$f_{k,l}:\Z_7^2\to \Z_7^2\quad\text{формулой}\quad f_{k,l}(x,y):=(2^kx,lx+y).$$
Проверьте, что

$\bullet$ таких преобразований ровно 21;

$\bullet$ они образуют группу;

$\bullet$ эта группа не является циклической (она даже не является абелевой).

{\it Другое указание.}
См. другое решение задачи  \ref{pr-cycle}-10.
Попробуйте сообразить, каким соотношениям должны удовлетворять перестановки $r$ и $s$,
чтобы множество $r^ks^l$, $k,l\in\Z$, образовывало бы нециклическую группу из 21 перестановки. А потом попробуйте придумать такие перестановки.

(b) Воспользуйтесь тем, что $2^5=33-1$.
Для любых целых неотрицательных $k,l$ определим преобразование $f_{k,l}:\Z_{11}^2\to \Z_{11}^2$ формулой $f_{k,l}(x,y):=(4^kx,lx+y)$.

(c) По теореме о первообразном корне существует элемент $a\in\Z_q$ порядка $p$.
Для любых целых неотрицательных $k,l$ определим преобразование $f_{k,l}:\Z_q^2\to \Z_q^2$ формулой $f_{k,l}(x,y):=(a^kx,lx+y)$.
Нетрудно проверить, что

$\bullet$ таких преобразований ровно $pq$;

$\bullet$ они образуют группу;

$\bullet$ эта группа не является циклической (она даже не является абелевой).

\bigskip
{\bf  \ref{pr-ordfh}.} (a) Поделите $n$ на $\ord f$ с остатком.

(b) Поделите $\ord h$ на $\ord_fh$ с остатком.

\smallskip
\textbf{\ref{pr-finite-ord}}. (d) Пусть $a$~--- элемент
порядка $3$. Выпишем все элементы группы. Теперь будем постепенно
зачеркивать их следующим образом: на каждом шаге выбираем
произвольным образом незачеркнутый элемент $x$ и зачеркиваем три элемента $x$, $xa$, $xa^2$.

При этом никакой элемент мы не зачеркнем больше одного раза.

(Действительно, предположим, что, например, зачеркиваемый элемент $xa$ уже был зачеркнут. Тогда либо
$xa=y$, либо $xa=ya$, либо $xa=ya^2$ для некоторого ранее
выбранного элемента $y$. Но тогда либо $x=ya^2$, либо $x=y$, либо
$x=ya$. Таким образом, элемент $x$ уже был зачеркнут. Противоречие. )

Значит, на каждом шаге зачеркивается ровно три новых элемента.

В конце будут зачеркнуты все элементы.
Значит, число элементов группы делится на $3$.

 (e) {\it Теорема Лагранжа.} Для $x\in G$ рассмотрим множество
$\{x,xf,xf^2,\dots,xf^{\ord f-1}\}$.
По определению порядка указанные элементы различны.
Значит, в этом множестве $\ord f$ элементов.
Если $xf^k=yf^l$, то $y=xf^{k-l}$.
Поэтому для разных $x$ эти множества либо не пересекаются, либо совпадают.
Значит, $|G|$ делится на $\ord f$.

\smallskip
{\bf  \ref{pr-comm}.} (a)
Обозначим через $p$ порядок неединичного элемента $f$.
Если $p=10$, то группа циклична.
Пусть теперь $p<10$.
По теореме Лагранжа $p\in\{5,2\}$.
Если есть элемент $g$ порядка~$10/p$, то $G=\{fg,(fg)^2,\dots,(fg)^{10}\}$.
Иначе есть элемент $g\not\in \{f,f^2,\dots,f^p\}$ порядка~$p$.
Тогда $\{f^kg^l\}_{k,l\in\mathbb{Z}}$ есть подгруппа из $p^2$ перестановок.
Противоречие с теоремой Лагранжа (задача \ref{pr-LagrTheorem}.b).

\smallskip
{\bf  \ref{pr-10}.} (a)
Если $f,g$ --- такие элементы, то $\{e,f,g,fg\}$ --- подгруппа из 4 перестановок.
Противоречие с теоремой Лагранжа (задача \ref{pr-LagrTheorem}.b).

\smallskip
{\bf  \ref{pr-abel}.} Все группы порядка $n$ коммутативны тогда и только тогда, когда
 в разложении $n=p_1^{k_1}\dots p_l^{k_l}$ числа $n$ на  простые сомножители
$k_i<3$ и $p_i$ не делит $p_j^{k_j}-1$.

Доказывается это примерно так же, как и для циклических групп, но в случае 1 надо воспользоваться тем, что любая конечная коммутативная группа раскладывается в прямую сумму циклических подгрупп.

\smallskip
\textbf{\ref{pr-orbits}}.
(a) Рассмотрим $25$ элементов $f^kg^l$ для $1\le k,l\le 5$.
Так как в группе всего $15$ элементов, то существуют $1\le k,l,m,n\le 5$ такие, что $(k,l)\ne (m,n)$ и
$f^kg^l=f^{m}g^n$.
Домножая на $f^{-m}$ слева и на $g^{-l}$ справа, получаем $f^{k-m}=g^{n-l}$.
Так как  $f$ и $g$~--- элементы порядка $5$, то множества $\{f,f^2,f^3,f^4\}$ и $\{g,g^2,g^3,g^4\}$ пересекаются.
%Из (a) получаем, что они совпадают.

(b) По (a) эти множества пересекаются.
Тогда существуют такие натуральные числа $1\le k,l\le 4$, что $f^k=g^l$.
Так как $\mathrm{\text{НОД}}(k,5)=1$, то существует такое целое число $m$, что $5\,|\,km-1$.
Тогда $f=f^{km}=(f^k)^m=(g^l)^m=g^{lm}$.
Отсюда следует, что эти множества совпадают.

(c) Следует из (b).

(d) Следуют из (c). (Пункт (d) не нужен для решения основной задачи.)
%, ибо порядки сопряженных перестановок равны (\ref{pr-fifteen}.b).

\smallskip
\textbf{\ref{pr-Lagr}}.
(a) Аналогично \ref{pr-orbits}.b.

(b) Иначе есть подгруппа из 9 элементов, что противоречит теореме Лагранжа.

\bigskip
{\bf  \ref{pr-1001_2}.} (a) Используйте задачу \ref{pr-ordfh}.ab.

(b) Следует из (a).

(c) Докажите, что $\ord(fh^{\ord f})$ делится на $q$ и на $\ord f$.

\smallskip
\textbf{\ref{br-type}.}  Перенумеруем элементы множества так, чтобы $f$ перешла в $g$.
Эта перенумерация задает требуемую перестановку $b$.

\smallskip
\textbf{\ref{pr-fifteen}.} (a,d) Утверждения доказываются индукцией по $k$.

\smallskip
{\bf  \ref{pr-alio}.} (d) Ввиду (c) $\widehat F=(s-1)15/s=15(1-1/s)$.

\smallskip
\textbf{\ref{pr-cen}}. (b) Для каждого $b\in G$ проведите стрелку от $f$ к $b^{-1}fb$.
Докажите, что если от $f$ к $g$ есть хотя бы одна стрелка, то общее количество стрелок от $f$ к $g$  равно
общему количеству стрелок от $f$ к $f$.

(c) Из \ref{pr-Lagr}.b следует, что с $f$ коммутируют только элементы $e, f, f^2$.
Теперь по (a) $c_G(f)=15/3=5$.

%(d) Иначе 14 делится на 4.

\smallskip
{\bf  \ref{pr-normalizer}.} (b) Ответ: $n(e)=n((ijk))=6$. 
%$S_3$
%и $n((ij))=?$.

(c) Аналогично \ref{pr-cen}.

\smallskip
\textbf{\ref{pr-pq}}. Аналогично разобранному случаю $p=3$, $q=5$.

\smallskip
\textbf{\ref{pr-case1}}. Аналогично \ref{pr-1001_2}.

\bigskip
{\bf  \ref{pr-main2}.} (d) Ввиду (c) $\widehat F=(|F|-|Z|)1001/|F|=1001(1-|Z|/|F|)$.

 \section{Доказательство из Книги}

\subsection{Формулировка основного результата}

{\it Группой} называется непустое семейство $G$ перестановок
%(т.е. преобразований)
некоторого множества, замкнутое относительно композиции и взятия обратной перестановки (т.е. если $f,g\in G$, то $f\circ g\in G$ и $f^{-1}\in G$).

\smallskip
{\bf Теорема (фольклор).} {\it
Для любой группы из $n$ перестановок найдется перестановка $g\in G$, для которой
$G= \{g,g^2,\dots,g^n\}$, тогда и только тогда, когда в разложении числа~$n$ на простые сомножители $n=p_1\dots p_t$

(*) все $p_i$ различны и

(**) $p_i$ не делит $p_j-1$ ни для каких $i$ и $j$.}

\smallskip
Заметим, что условие `(*) и (**)' равносильно взаимной простоте $n$ с $\phi(n)$.
Здесь $\phi(n)$  --- количество целых чисел от 1 до $n$, взаимно простых с $n$ (функция Эйлера).

Интересные частные случаи приведены в задачах \ref{pr-prime-cycle}, \ref{pr-pq} и \ref{pr-constr}.c.
Как его придумать, видно из \S2.
%добавления, ср. [BKSS].

Если в конечной группе $G$ найдется перестановка $g$, из всех возможных степеней которой состоит $G$
(т.е. $G=\{g,g^2,\dots,g^n,\dots\}$), то эта группа называется {\it циклической\/}.

%\newpage
\subsection{Доказательство части <<только тогда>>}

Обозначим $\Z/k:=\{(1,2,\dots,k)^i\ |\ i=1,2,\dots,k\}$.

Несвязное объединение множеств $M$ и $N$ определяется как $M\sqcup N:=M\times\{0\}\cup N\times\{1\}$.
Если $g$ и $h$ --- преобразования множеств $M$ и $N$, то  преобразование
$g\sqcup h$ несвязного объединения $M\sqcup N$ определяется формулой
$(g\sqcup h)(x):=\begin{cases}g(x) & x\in M\times\{0\}\subset M\sqcup N\\
h(x) & x\in N\times\{1\}\subset M\sqcup N\end{cases}$.
%декартова произведения $M\times N$
%Обозначим $G\oplus H=\{g\times h\ |\ g\in G,\ h\in H\}$.
Для групп $G$ и $H$, состоящих из преобразований множеств $M$ и $N$, определим
$$G\times H:=\{g\sqcup h: M\sqcup N\to  M\sqcup N\ |\ g\in G,\ h\in H\}.$$
Если нарушается вышеприведенное условие (*), например, $p_1=p_2=p$, то в качестве
нециклической группы из $n$ элементов можно взять группу $\Z/p\times\Z/\frac np$.
\footnote{Т.е. группу
$\left\{(1,2,\dots,p)^i(p+1,p+2,\dots,p+\tfrac np)^k\ |\ i=1,\dots,p,\ k=1,\dots,\tfrac np\right\}.$}

Обозначим через $\Z_k$ множество вычетов по модулю $k$
с операциями суммы и произведения, известными из теории чисел.
(Не путайте с группой $\Z/k$, определенной выше!)
Если нарушается вышеприведенное условие (**), например, $p_1$ делит $p_2-1$, то
{\it существует вычет $a\in\Z_{p_2}$, для которого степени $a,a^2,\dots,a^{p_1}=1$ различны}.
Это следует из теоремы о первообразном корне (сформулированной в \S1).
% (если Вы ее не знаете, примите указанное следствие на веру).

Обозначим через $G_{p_1,p_2}$ группу преобразований $f_{k,l}:\Z_{p_2}^2\to \Z_{p_2}^2$, заданных  формулой $f_{k,l}(x,y):=(a^kx,lx+y)$ для $k\in\Z_{p_1}$ и $l\in\Z_{p_2}$.
\footnote{
Научно говоря,
$G_{p_1,p_2}=\left\{
\begin{pmatrix}
a^k&l\cr
0&1\cr
\end{pmatrix}\in\Z_{p_2}^{2\times2}
\;\Big|\;
k\in\Z_{p_1},\ l\in\Z_{p_2}\;
\right\}.$}
%коммутатив
Тогда в качестве нециклической (даже некоммутативной) группы из $n$ элементов можно взять группу
$G_{p_1,p_2}\times\Z/\frac n{p_1p_2}$.
\footnote{Т.е. группу
$\left\{f\circ(1,2,\dots,\frac n{p_1p_2})^j\ |\ f\in G_{p_1,p_2},\ j=1,2,\dots,\frac n{p_1p_2}\right\}.$}
QED

 %\newpage
\subsection{Доказательство части <<тогда>>}

Через $|X|$ обозначается число элементов в множестве $X$.
Обозначим данную группу через $G$.
Используем индукцию по числу простых сомножителей в $n=|G|$.
Если сомножитель один, то часть <<тогда>> вытекает из следующего частного случая теоремы Лагранжа.

{\bf Порядком} $\ord a$ элемента $a$ группы с единичным элементом $e$ называется наименьшее
целое положительное $n$, для которого $a^n=e$.
Если группа конечна, то ясно, что такое $n$ существует.

\smallskip
{\bf Теорема Лагранжа (частный случай).} {\it Число элементов конечной группы делится на порядок любого ее элемента.}

\smallskip
{\it Доказательство.}
\footnote{Более подробно это доказательство (и его обобщение, см. ниже) изложено в [KS, стр. 64].}
Обозначим данную группу через $G$.
Для любого $x\in G$ рассмотрим множество
%\linebreak
$\{x,xf,xf^2,\dots,xf^{\ord f-1}\}$.
Из определения порядка вытекает, что указанные элементы различны.
Значит, в этом множестве $\ord f$ элементов.
Если $xf^k=yf^l$, то $y=xf^{k-l}$.
Поэтому для разных $x$ эти множества либо не пересекаются, либо совпадают.
Значит, $|G|$ делится на $\ord f$.
QED

\def\gp#1{\left\langle#1\right\rangle}

\smallskip
Пусть теперь простых сомножителей в $n=|G|$ больше одного.
Нам понадобится следующая общая версия теоремы Лагранжа.

{\bf Подгруппой} группы называется подмножество этой группы, которое само по себе является группой.

\smallskip
{\bf Теорема Лагранжа.} {\it Число элементов конечной группы делится на число элементов любой ее подгруппы.}

\smallskip
{\it Доказательство.} Обозначим данную группу через $G$, а ее подгруппу через $\{h_1,\dots,h_m\}$.
Для любого $x\in G$ рассмотрим множество
%\linebreak
$\{xh_1,xh_2,\dots,xh_m\}$.
В этом множестве $m$ элементов.
Если $xh_k=yh_l$, то $y=xh_kh_l^{-1}$.
Поэтому для разных $x$ эти множества либо не пересекаются, либо совпадают.
Значит, $|G|$ делится на $m$.
QED

%{\bf Vyalyi: неясно, зачем давать два этих доказательства. Их трудность
%(при чтении) примерно одинакова и точно такая же, как при чтении
%учебника по теории групп. Из второе следует первое (и ниже все равно
%вводится группа, порожденная элементом).}

\smallskip
{\bf Максимальной подгруппой} назовем максимальную по включению подгруппу, не совпадающую со всей группой и содержащую более одного элемента.
По предположению индукции и теореме Лагранжа {\it каждая максимальная подгруппа является циклической.}

Для элемента $f$ группы $G$ обозначим через $\gp f\subset G$ множество всех его степеней (в т.ч. нулевой и отрицательных).
Элемент $f$ называется {\bf порождающим} для (циклической) подгруппы $\gp f$.

Предположим противное, т.е. что группа $G$ не является циклической.
Тогда {\it каждый элемент $f$ содержится в некоторой максимальной подгруппе} (в максимальной по включению подгруппе, содержащей $\gp f$).

Элементы $f$ и $g$ группы $G$ называются {\bf сопряженными} в $G$, если $g=b^{-1}fb$ для некоторого $b\in G$.

%{\bf Vyalyi: в этом месте даже для МП я бы вставил упоминание о действии группы сопряжениями.
%При выбранном подходе это мучительно сложно сделать, так как группа у вас - это уже какие-то преобразования и
%действовать по-другому она не может.
%Но суть вашего доказательство именно в действии группы сопряжениями.
%АС: Не совсем, нужны лишь его орбиты.
%Но мы готовы вставить текст, который кто-то предложит.}

\smallskip
{\bf Первый случай:}
{\it порождающий элемент $f$ некоторой максимальной подгруппы сопряжен только с некоторыми своими степенями.}
Возьмем $h\in G-\gp f$.
Обозначим через $q$ наименьшее из целых положительных $n$, для которых $h^n\in \gp f$.
Такое $q$ существует, поскольку $h^{\ord h}\in \gp f$.

%Тогда $h^{-1}f^sh=f^{ks}$ для любого $s\in\Z$.

\smallskip
{\it Доказательство того, что $|G|$ делится на $q$.}
Поделим с остатком $\ord h$ на $q$: $\ord h=qt+r$.
Тогда $h^r=h^{\ord h-qt}\in \gp f$ и $0\le r<q$.
Поэтому ввиду минимальности $q$ имеем $r=0$.
Т.е. $\ord h$ делится на $q$.
Значит, по теореме Лагранжа $|G|$ делится на $q$.
\qed

\smallskip
{\it Доказательство того, что $hf=fh$.}
По условию первого случая $h^{-1}fh=f^k$ для некоторого $k\in\Z$.
Значит, $f=h^{-q}fh^q=f^{k^q}$
(последнее равенство верно для произвольного $q$ и доказывается индукцией по $q$).
Поэтому $k^q\equiv1\mod\ord f$.
Следовательно, $k$ и $\ord f$ взаимно просты.
По теореме Лагранжа и условиям (*) и (**), $|G|$ и $\varphi(\ord f)$ взаимно просты.
Так как $|G|$ делится на $q$, то $q$ и $\varphi(\ord f)$ взаимно просты.
Следовательно, найдутся целые $x$ и $y$, для которых $qx+\varphi(\ord f)y=1$.
Значит, $k\equiv k^{qx+\varphi(\ord f)y}\equiv1\mod\ord f$.
Поэтому $fh=hf$.
QED

%$\varphi(\ord f)$ делится на $q$.?
%Значит, их общий делитель $q$ равен единице.

%\footnote{Вот немного другое доказательство утверждения $fh=hf$, предложенное М.Н. Вялым.
%Обозначим через $r$ наименьшее из целых положительных $n$, для которых $h^{-n}fh^n=f$.
%Такое $r$ существует, поскольку $h^{-n}fh^n=f$ для $n=\ord h$.
%Дальнейшие рассуждения буквально повторяют вышеизложенное доказательство с заменой $r$ на $q$. }

%По условию (*) и теореме Лагранжа $\ord f$ является произведением $p_1\dots p_s$ различных простых.
%Тогда $k^q\equiv1\mod p_i$ для любого $i=1,2,\dots,s$.
%Так как $|G|$ делится на $\ord h$ и $\ord h$ делится на $q$, то $|G|$ делится на $q$.
%Поэтому и по условию (*) $q$ является произведением различных простых.
%По условию (**) ни одно из этих простых $p_j$ не делит никакое $p_i-1$.
%Следовательно, $q$ взаимно просто с каждым $p_i-1$.
%Поэтому существуют целые $x=x_i$ и $y=y_i$, для которых $qx+(p_i-1)y=1$.
%Значит, $k\equiv k^{qx+(p_i-1)y}\equiv1\mod p_i$ для любого $i=1,2,\dots,s$.

\smallskip
{\it Завершение разбора первого случая.}
Так как $fh=hf$, то в $G$ есть подгруппа
$$\{f^ih^j\ |\ 1\le i\le \ord f,\ 1\le j\le q\}$$
из $q\ord f$ элементов.
Значит, по условию (*) и теореме Лагранжа $\ord f$ и $q$ взаимно просты.
Обозначим $g:=fh^{\ord f}$.
Так как $g^j=f^jh^{j\ord f}$ для любого $j$, то $\ord g$ делится на $q$ и на $\ord f$.
Поэтому $\ord g=q\ord f$.
Так как подгруппа $\gp f$ максимальна, то $\gp g=G$.
Значит, $G$ циклическая. Противоречие.
QED

%Vyalyi (к следующему): фактически вы смотрите на нормализатор максимальной группы и
%две возможности для него. Первый случай: нормализатор некоторой
%максимальной группы совпадает с
%G. Почему бы так и не сказать? Дальнейшее ваше рассуждение на обычном
%языке звучит так: ``если максимальная группа нормальна, то она лежит
%в центре (действие сопряжением любым элементом группы G имеет
%порядок, который делит и порядок G, и $\phi$(порядка М), т.е. этот
%порядок 1), а тогда и вся группа циклическая.
%AS: последнего я не понял. Случай коммутативной группы не очевиден и использует (*).

\smallskip
%\newpage
{\bf Второй случай:}
{\it порождающий элемент любой максимальной подгруппы сопряжен не только со своими степенями.}
\footnote{Вот план немного другого разбора второго случая, предложенный М.Н. Вялым.
Сначала вводим множество $N(F)$ из доказательства утверждения (3) ниже.
Доказываем, что  $N(F)=F$, см. там же.
Тогда (1) очевидно, так $Z(G)\subset N(F)$.
Далее делаем то же, что и в приводимом разборе. }

{\bf Произведением двух подмножеств $X$ и $Y$} группы $G$ называют
множество всевозможных произведений $xy$, где $x\in X$ и $y\in Y$.
Если одно из этих подмножеств состоит только из одного элемента,
например, $Y=\{y\}$, то для краткости пишут $Xy$ вместо
$X\{y\}$.

{\bf Центром} группы $G$ называют
$$Z=Z(G):=\{a\in G\ :\ ga=ag\text{ для любого }g\in G\},$$
т.е. множество тех элементов, которые коммутируют со всеми.

%:=\{fz\in G\ :\ f\in F,\ z\in Z\}$

\smallskip
{\it (1) Любая максимальная подгруппа содержит центр. }

\smallskip
{\it Доказательство утверждения (1).}
Пусть, напротив, максимальная подгруппа $F$ не содержит $Z$.
Тогда $FZ$ --- б\'ольшая подгруппа, чем $F$.
Ввиду максимальности подгруппы $F$ имеем $FZ=G$.
Значит, $G$ коммутативна, т.е. $G=Z(G)$.
Поэтому порождающий элемент любой максимальной подгруппы сопряжен только с собой.
Это противоречит условию второго случая. QED

\smallskip
{\it (2) Пересечение двух максимальных подгрупп равно центру.}

\smallskip
{\it Доказательство утверждения (2).}
%коммутатив
Неединичный элемент в пересечении коммутирует с элементами обоих подгрупп.
Значит, он коммутирует с любым произведением нескольких сомножителей, каждый из которых
лежит в одной из наших подгрупп.
Множество таких произведений является подгруппой, содержащей обе максимальные подгруппы.
В силу максимальности наших подгрупп эта подгруппа совпадает со всей группой.
Значит, пересечение содержится в центре.

Из (1) вытекает обратное включение.
QED

\smallskip
{\it Завершение разбора второго случая: подсчет.}
Напомним, что любой элемент группы содержится в некоторой максимальной подгруппе.
Ввиду (2) для любого {\it нецентрального} элемента такая подгруппа единственна.
Поэтому вся группа разбивается на центр и несвязное объединение дополнений максимальных подгрупп до центра.
Количества элементов в таких дополнениях одинаковы для сопряженных максимальных подгрупп.
(Подмножества $F,F'\subset G$ называются {\it сопряженными}, если $g^{-1}Fg=F'$ для некоторого $g\in G$. )
%всех нецентральных элементов, сопряженных некоторому элементу максимальной подгруппы $F$,
Обозначим через $\widehat F$ общее количество нецентральных элементов во всех максимальных подгруппах,
сопряженных с максимальной подгруппой $F$.
Обозначим через $F_1,\dots,F_s$ наибольший набор попарно несопряженных максимальных подгрупп.
%Тогда будем рассматривать классы сопряжённости максимальных подгрупп. Пусть их $k$.
%Тогда пусть в каждой подгруппе $i$-того класса сопряжённости $F_i$ элементов.
Тогда
$$|G|=|Z|+\sum_{i=1}^s\widehat{F_i}.$$
Ввиду левого неравенства в следующем утверждении число слагаемых не превосходит единицы, а ввиду правого
--- одного слагаемого тоже быть не может.
QED

\smallskip
{\it (3) $|G|/2\le \widehat F<|G|-|Z|$.}

\smallskip
{\it Доказательство утверждения (3).}
Подгруппа, сопряженная к максимальной, также максимальна.

(Действительно, если $g^{-1}Fg\subset F'\subset G$, то $F\subset gF'g^{-1}\subset G$.)

Рассмотрим множество
$$N(F):=\{a\in G\ :\ Fa=aF\}.$$
%Пусть $F=\gp f$.
%Рассмотрим множество
%$$N(f):=\{a\in G\ :\ fa=af^k\text{ для некоторого }k\}.$$
Тогда число различных подгрупп, сопряженных с $F$ (включая $F$), равно $|G|/|N(F)|$.
\footnote{Это вытекает из теоремы о длине орбиты для действия группы на себе сопряжениями.
Те, кому это доказательство непонятно, могут прочитать следующий абзац.}

(Действительно, сопряжение каждым элементом группы $G$ переводит подгруппу $F$ в одну из сопряженных подгрупп.
Если сопряжение двумя разными элементами $u,v$ группы $G$ переводит подгруппу $F$ в одну и ту же подгруппу, т.е.
$u^{-1}Fu=v^{-1}Fv$, то $(uv^{-1})^{-1}Fuv^{-1}=F$.
Это означает, что $uv^{-1}\in N(F)$ или, что то же самое, $u\in N(F)v$.
Обратно, условие $u\in N(F)v$ влечет $u^{-1}Fu=v^{-1}Fv$.
Ясно, что $|N(F)v|=|N(F)|$.
Поэтому число элементов в $G$, сопряжение с которыми переводит подгруппу $F$ в данную фиксированную сопряженную подгруппу, равно $|N(F)|$.
Значит, число подгрупп, сопряжeнных к $F$, ровно в $|N(F)|$ раз меньше, чем элементов группы $G$. )

Имеем $N(F)=F$.

(Действительно, $N(F)$ является подгруппой.
По условию второго случая $N(F)\ne G$.
Так как $N(F)\supset F$, то в силу максимальности $N(F)=F$.)

Ввиду трех доказанных утверждений $\widehat F=(|F|-|Z|)\dfrac{|G|}{|F|}=|G|\left(1-\dfrac{|Z|}{|F|}\right)$.

Так как $|G|>|F|$, то $\widehat F<|G|-|Z|$.

По (1) центр является подгруппой в $F$.
Значит, по теореме Лагранжа $|Z|$ делит $|F|$.
По условию второго случая $Z\ne F$.
Следовательно, $|Z| \le |F|/2$.
Поэтому $\widehat F\ge|G|/2$.
QED

\bigskip
%\newpage
\centerline{\bf Литература}

[A] В.И. Арнольд, Обыкновенные дифференциальные уравнения, М, Наука, 1984.

[Al] В.Б. Алексеев, Теорема Абеля. М: Наука, 1976.

[B] Ken Brown,
Mathematics 4340, When are all groups of order $n$ cyclic?
Cornell University, March 2009,
\url{http://www.cornell.edu/~kbrown/4340/cyclic_only_orders.pdf}

[BKKSS] Д. Баранов, А. Клячко, K. Кохась, А. Скопенков и М. Скопенков,
Когда любая группа из $n$ элементов циклическая?
\url{http://olympiads.mccme.ru/lktg/2011/6/index.htm}

[GDI] А.А. Глибичук, А.Б. Дайняк, Д.Г. Ильинский, А.Б. Купавский, А.М. Райгородский, А.Б. Скопенков, А.А. Чернов,
Элементы дискретной математики в задачах, Изд-во МЦНМО, 2015.

[K] Ф. Клейн, Элементарная математика с точки зрения высшей.

[KS] Л.А. Калужнин и В.И. Сущанский, Преобразования и перестановки, М.: Наука, 1985.

[P1] А. Пуанкаре, О науке, М.: Наука, 1990.

[P2] А.Пуанкаре `Наука и метод',

[S08] А. Скопенков, Еще несколько доказательств из Книги: разрешимость
и неразрешимость уравнений в радикалах, \url{http://arxiv.org/abs/0804.4357}

[S15] A. Skopenkov, A short elementary proof of the Ruffini-Abel Theorem
\linebreak
\url{http://arxiv.org/abs/1508.03317}

[V] Виноградов И.М., Основы теории чисел. М.; Ижевск: НИЦ <<Регулярная и хаотическая динамика>>, 2003.

[Z] Математика в задачах. Сборник материалов московских выездных математических школ.
Под редакцией А. Заславского, Д. Пермякова, А. Скопенкова, М. Скопенкова и
А. Шаповалова. Москва, МЦНМО, 2009.
\url{http://www.mccme.ru/circles/oim/mvz.pdf}

\end{document}